\documentclass[12pt]{article}
\usepackage{amsfonts,amssymb,amsmath}
\usepackage{graphicx}

\usepackage{color} 
   \definecolor{cites}{rgb}{0.50 , 0.00 , 0.00}  
   \definecolor{urls} {rgb}{0.00 , 0.00 , 0.50}  
   \definecolor{links}{rgb}{0.00 , 0.00 , 0.50}   
\usepackage[
      colorlinks=true,   
      citecolor=cites,   
      urlcolor=urls,     
      linkcolor=links,   
      pdftitle={The Main Diagonal of a Permutation Matrix},  
      pdfauthor={Marko Lindner, Gilbert Strang}, 
      pdfpagemode=UseOutlines, 
      pdfstartview=FitH,       
      bookmarksopen=false      
   ]{hyperref}

\newcommand\eps\varepsilon
\newcommand\ph\varphi

\newcommand\codim{{\rm codim}}

\newcommand\ind{{\rm index}}

\newcommand\rank{{\rm rank}}

\newcommand\Z{{\mathbb Z}}
\newcommand\N{{\mathbb N}}
\newcommand\I{{\mathbb I}}

\newcommand\obullet{{\bigcirc\hspace{-1.75ex}\bullet}}

\newtheorem{bigtheorem}{Theorem}  
\newtheorem{theorem}{Theorem}[section]
\newtheorem{lemma}[theorem]{Lemma}
\newtheorem{corollary}[theorem]{Corollary}

\newenvironment{remark}
 {\par\noindent\refstepcounter{theorem}{\bf Remark \thetheorem}\ }
 {\raisebox{1mm}{\framebox{}}\pagebreak[2]}

\newenvironment{example}
 {\par\noindent\refstepcounter{theorem}{\bf Example \thetheorem}\ }
 {\raisebox{1mm}{\framebox{}}\pagebreak[2]}

\newcommand\proofend{{\rule{2mm}{2mm}}}

\newenvironment{proof}
 {\par\noindent{\bf Proof.}}
 {\proofend\pagebreak[2]}

\numberwithin{figure}{section}  

\begin{document}
\title{The Main Diagonal of a Permutation Matrix}
\author{
{\sc Marko Lindner}\quad and\quad {\sc Gilbert Strang}}
\date{~}
\maketitle


\renewcommand{\baselinestretch}{1.0}
{
\footnotesize\begin{center}\sc Abstract\end{center}

\noindent By counting 1's in the ``right half'' of $2w$ consecutive rows, we locate the main diagonal of any doubly infinite permutation matrix with bandwidth $w$. Then the matrix can be correctly centered and factored into block-diagonal permutation matrices.

Part \ref{part:2} of the paper discusses the same questions for the much larger class of band-dominated matrices. The main diagonal is determined by the Fredholm index of a singly infinite submatrix. Thus the main diagonal is determined ``at infinity'' in general, but from only $2w$ rows for banded permutations.

\noindent
{\it Mathematics subject classification (2010):} 15A23; 47A53, 47B36.\\
{\it Keywords and phrases:} banded matrix, permutation, infinite matrix, main diagonal, factorization.
}
\section{Introduction}
This paper is about banded doubly infinite permutation matrices. These matrices represent permutations of the integers $\Z$, in which no integer moves more than $w$ places. The banded matrix $P$ has zero entries $p_{ij}$ whenever $|i-j|>w$, and it has exactly one entry $p_{ij}=1$ in each row and column. So $P$ has $2w+1$ diagonals that are possibly nonzero, and exactly one of these deserves to be called the {\sl main diagonal}. The first objective of this paper is to find the main diagonal of $P$.

The remarkable fact is that the correct position of the main diagonal is determined by {\sl any $2w$ consecutive rows of $P$}. This will be Theorem 1. Normally that position can only be determined ``at infinity''.

For a much larger class of biinfinite matrices $A$, we form the singly infinite submatrix $A_{+}$ from rows $i>0$ and columns $j>0$ of $A$. The {\sl plus-index} of $A$ is the usual Fredholm index of $A_{+}$, computed from the dimension of its nullspace and the codimension of its range:
\begin{equation} \label{eq:ind+}
\ind_{+}(A)\ =\ \ind(A_{+})\ =\ \dim N(A_{+})\,-\,\codim\,R(A_{+}).
\end{equation}
For any Fredholm operator, including all permutations $A=P$, those two numbers are finite. Then the main diagonal of $A$ is determined by $\kappa=\ind_{+}(A)$. It is found $\kappa$ diagonals above (or $-\kappa$ diagonals below) the zeroth diagonal of $A$. So our problem is to compute that index.

\begin{example} \label{ex:shift1}
The doubly infinite forward shift $S$ has nonzero entries $S_{i,i-1}=1$ for $i\in\Z$. Its singly infinite submatrix $S_{+}$ is lower triangular, with those ones along the first subdiagonal. The nullspace of $S_{+}$ has dimension $=0$ (independent columns of $S_{+}$) but the range has codimension $=1$ (it consists of all singly infinite vectors with a zero in the first position). Thus $\kappa=\ind_{+}(S)=0-1$ and the main diagonal of $S$ is correctly located: It is one diagonal below the zeroth diagonal of $S$ (and it contains the ones).

Since $S$ is a permutation of $\Z$ with bandwidth $w=1$, Theorem 1 says that $\ind_{+}(S)$ can be found from any two consecutive rows $i$ and $i+1$ of $S$: ``Count the ones in columns $i+1,\,i+2,\,\dots$ and subtract $w$.'' The result $0-1$ agrees with $\ind_{+}(S)$.
\end{example}

We can directly state Theorems 1--3, for banded permutations $P$ of $\Z$. Two proofs are given for Theorem 1. Then Part \ref{part:2} of the paper discusses the much larger class of doubly infinite band-dominated Fredholm operators $A$. The plus-index of $A$ (and thus its main diagonal) is well defined but not so easily computable.

\begin{bigtheorem} 
If $P$ has bandwidth $w$, its plus-index is determined by rows $1,\dots,2w$ and columns $j>w$. Subtract $w$ from the number $n$ of ones in this submatrix. The result $\kappa=n-w$ is the plus-index of $P$, and the main diagonal of $P$ is $\kappa$ diagonals above the zeroth diagonal.

The same result $\kappa$ comes from rows $j^{*}-w$ to $j^{*}+w-1$ and columns $\geq j^{*}$, for any $j^{*}\in\Z$. The submatrix is shown in~\eqref{eq:submatrix} below.
\end{bigtheorem}

Equivalently, $P_c=S^\kappa P$ is a {\sl centered} permutation of $\Z$, by which we mean that its zeroth diagonal is the main diagonal. Our first proof will show how the submatrix with $2w$ rows was discovered. The second proof goes directly to the index $\kappa$.


\begin{bigtheorem} 
Every centered permutation $P_c$ of $\Z$ with bandwidth $w$ can be factored into a product of two block-diagonal permutations:
\begin{equation} \label{eq:fact2}
P_c\ =\ BC\ =\ \left(
\begin{array}{cccccccc}
\cdot\\[-0.7ex]
\cline{2-4}
& \multicolumn{3}{|c|}{~}\\[-1.5ex]
& \multicolumn{3}{|c|}{B_{0}}\\[-1.5ex]
& \multicolumn{3}{|c|}{~}\\
\cline{2-7}
& &&& \multicolumn{3}{|c|}{~}\\[-1.5ex]
& &&& \multicolumn{3}{|c|}{B_{1}}\\[-1.5ex]
& &&& \multicolumn{3}{|c|}{~}\\
\cline{5-7}
& &&& &&& \ddots~\\[-1.5ex]
~
\end{array}
\right)
 \left(
\begin{array}{cccccccc}
~\\[-1.5ex]
~\ddots\\
\cline{2-4}
& \multicolumn{3}{|c|}{~}\\[-1.5ex]
& \multicolumn{3}{|c|}{C_{0}}\\[-1.5ex]
& \multicolumn{3}{|c|}{~}\\
\cline{2-7}
& &&& \multicolumn{3}{|c|}{~}\\[-1.5ex]
& &&& \multicolumn{3}{|c|}{C_{1}}\\[-1.5ex]
& &&& \multicolumn{3}{|c|}{~}\\
\cline{5-7}
& &&& &&& \cdot
\end{array}
\right)
\end{equation}
All blocks have size $2w$ and each $C_{i}$ is ``offset'' between $B_{i}$ and $B_{i+1}$ (shifted by $w$ rows and columns, as shown).
\end{bigtheorem}

\begin{bigtheorem} 
The centered permutation $P_{c}$ can be further factored into $N<2w$ block-diagonal permutations of bandwidth $1$:
\begin{equation} \label{eq:factN}
P_{c}\ =\ F_{1}\,F_{2}\,\cdots\,F_{N}
\end{equation}
\end{bigtheorem}

Each factor $F_{i}$ has block size $1$ or $2$. Thus $F_i$ exchanges disjoint pairs of neighboring rows. By Theorems 1--3 the original $P$ is factored into
\[
P\ =\ S^{-\kappa}P_{c}\ =\ S^{-\kappa}BC\ =\ S^{-\kappa}F_{1}F_{2}\cdots F_{N}.
\]

\section{The index and the plus-index} 

We say that an infinite matrix $A=(a_{ij})_{i,j\in\Z}$ is {\sl invertible} if the linear operator that it represents via matrix-vector multiplication is invertible as an operator from $\ell^2(\Z)$ to $\ell^2(\Z)$.
A bounded linear operator $A$ from $\ell^{2}$ into $\ell^{2}$ is invertible iff it is both injective (its nullspace $N(A)=\{x\in \ell^{2}:Ax=0\}$ consists of $0$ only) and surjective (its range $R(A)=\{Ax:x\in \ell^{2}\}$ is all of $\ell^{2}$). Deviation from both properties is measured in terms of the two integers
\begin{equation} \label{eq:alpha_beta}
\alpha\ =\ \dim N(A)\qquad\textrm{and}\qquad\beta\ =\ \codim\,R(A).
\end{equation}
$A$ is a {\sl Fredholm operator} if both $\alpha$ and $\beta$ are finite. Then the Fredholm index (or just the {\sl index}) of $A$ is the difference
\begin{equation} \label{eq:index}
\ind(A)\ =\ \alpha\, -\, \beta.
\end{equation}
Unlike the separate numbers $\alpha$ and $\beta$, their difference $\alpha-\beta$ has these important properties:
\begin{itemize}
\item the index is invariant under compact perturbations $A+K$;
\item the index of a product obeys the remarkable formula $\ind(AB)\,=\,\ind(A) +\ind(B)$;
\item the index is continuous with respect to the operator norm of $A$ (and therefore locally constant); and
\item all finite square matrices have index zero.
\end{itemize}
All Fredholm matrices of index zero, including all finite square matrices, obey the {\sl Fredholm alternative}: Either
$A$ is injective {and} surjective ($\alpha=0$, $\beta=0$)
or
$A$ is not injective and not surjective ($\alpha\ne 0$, $\beta\ne 0$).
This is the set of all operators $A=C+K$ with $C$ invertible and $K$ compact. See e.g. \cite{Davies2007:Book,GohGoldKash} for a nice introduction to Fredholm operators. The property of being invertible on $\ell^p(\Z)$ and the index itself are independent of $p\in [1,\infty]$ for a certain class of infinite matrices (with uniformly bounded entries and summable off-diagonal decay), see \cite{Kurbatov,Li:Wiener}.

Our interest in the index (more precisely, the plus-index) of a biinfinite matrix originates from the following natural question:
\begin{center}
{\sl Which diagonal is the main diagonal of a biinfinite matrix?}
\end{center}
For a symmetric matrix, the zeroth diagonal is the main diagonal. For a Toeplitz matrix (a polynomial in the shift $S$), a zero winding number is the key. A wider class of structured matrices was analyzed by de Boor \cite{deBoor}. For Fredholm operators in general, Israel Gohberg's diplomatic answer to this question was that every diagonal has the right to be the main diagonal \cite[p. 24]{Bo96}. But there are concrete problems waiting for a concrete answer. Here are two such problems:

\begin{enumerate}
\item For finite and semiinfinite matrices, the inverse of a lower triangular matrix is again lower triangular. This may fail if $A$ is biinfinite (essentially because the `wrong' diagonal is mistaken for the main diagonal). The inverse of the lower triangular shift $S$ is the upper triangular backward shift $S^{\top}$. Shifting $S$ one row up (treating the ones as the main diagonal) resolves this conflict.

Here is a slightly more sophisticated example: $A=S-\frac12 S^2$ is lower triangular with its inverse neither lower nor upper triangular. Shift up by one row: $B=S^{-1}A=I-\frac 12 S$ is lower triangular with $B^{-1}=I+\frac 12 S+\frac14 S^2+\cdots$ also lower triangular (as we want). Shift up one more row: $C=S^{-2}A=S^{-1}-\frac12 I$ is now upper triangular with $C^{-1}=B^{-1}S$ lower triangular. In this example the diagonal consisting of all ones should be the main diagonal, and $B$ is centered (even though triangular) so that its inverse is triangular of the same kind.
\item The numerical solution of biinfinite systems $Ax=b$ can approximate $A^{-1}$ by the inverses of finite square submatrices ({\sl finite sections}) of $A$. Their upper left and lower right corners lie on the main diagonal of $A$. But again, which diagonal should that be? For $A=2S^2+\frac12 S-I+2S^{-1}$, one can show \cite{BGr5,BSi2} that it has to be the diagonal that carries all the $\frac 12$'s. Finite sections that are centered along one of the other nonzero diagonals are invertible (for sufficiently large sizes) but their inverses do not converge (they blow up).
\end{enumerate}
Problem $1$ and $2$ are not unrelated: If $A$ is lower triangular and the inverses of its finite sections (which are also lower triangular) approximate $A^{-1}$, then $A^{-1}$ will be lower triangular, too.
A similar argument can be used to transfer Asplund's theorem \cite{Asplund,Strang:Asplund} to infinite matrices, explaining the relations between ranks of submatrices of $A$ and $A^{-1}$. We demonstrate this in Part \ref{part:2} of this paper.

The answer to both problems is the same:
\textsl{Shift $A$ down by $\kappa$ rows}, or equivalently by $\kappa$ diagonals,
where $\kappa$ is the {\sl plus-index} of $A$ \cite{RaRoRoe}:
\begin{equation} \label{eq:plusindex}
\kappa\, =\, \ind_{+}(A)\, =\, \ind(A_+)\qquad\textrm{and}\qquad
A_+\, =\, \left(\begin{array}{ccc}
a_{11}&a_{12}&\cdots\\
a_{21}&a_{22}&\cdots\\
\vdots&\vdots&\ddots
\end{array}\right)
\end{equation}
$A_+$ is a semiinfinite submatrix of $A$. This shifting process is called {\sl index cancellation} \cite{GohbergFeldman,HeinigHellinger}.

The centered matrix $A_{c}=S^\kappa A$ has $\ind_{+}(A_{c})=0$ (see \eqref{eq:IC} below). This is necessary for a triangular matrix to have a triangular inverse of the same kind (Corollary \ref{cor:triangular}). It is also necessary for the convergence of the inverses of the finite sections to $A^{-1}$ (see \cite{Li:FSMsubs,SeidelSilbermann3}).

How is this plus-index $\kappa$ computed? In \cite{RaRoRoe} it is shown that $\kappa$ is invariant under passing to a ``limit operator'' \cite{RaRoSiBook,Li:Book,CWLi2008:Memoir} of $A$ at $+\infty$, which often simplifies the computation. For the plus-index of a permutation $P$, the limit operator is not needed and the new formula \eqref{eq:indfinal} is as simple as possible.

\part{Permutation matrices} \label{part:1}
\section{The plus-index of a biinfinite permutation matrix} 
Now we come to the problem of computing the plus-index of a biinfinite permutation matrix. So let $\pi:\Z\to\Z$ be a permutation (a bijection) of the integers and put $P=(p_{ij})_{i,j\in\Z}$ with
\[
p_{ij}\ =\ \delta_{\pi(i),j}\ =\ \left\{\begin{array}{cl}
1&\textrm{if } j=\pi(i),\\
0&\textrm{if } j\ne \pi(i).
\end{array}\right.
\]

The matrix $P$ is banded iff the maximal displacement $w$ of any integer via $\pi$ is {\sl finite}:
\begin{equation} \label{eq:w}
w\ :=\ \sup_{i\in\Z}|i-\pi(i)|\quad\textrm{is the bandwidth of $P$.}
\end{equation}

Every permutation matrix $P$ is invertible. So the submatrix $P_+$ is always Fredholm (see Lemma \ref{lem:A+-} below) and it makes sense to ask for its index: the plus-index of $P$.


It seems a bit arbitrary to define the plus-index of a biinfinite matrix $A$ based on the submatrix $A_{+}$ that starts at the particular entry $a_{11}$. Lemma \ref{lem:Ak} shows that for banded permutations, the submatrix $A_k$ starting at $a_{kk}$ gives the same plus-index for every  $k\in\Z$:
\begin{equation} \label{eq:defAk}
A_k\ :=\ (a_{ij})_{i,j=k}^\infty\ =\
\left(\begin{array}{ccc}
a_{k,k}&a_{k,k+1}&\cdots\\
a_{k+1,k}&a_{k+1,k+1}&\cdots\\
\vdots&\vdots&\ddots
\end{array}\right)
\end{equation}

\begin{lemma} \label{lem:Ak}
If $P$ is a banded biinfinite permutation matrix and $k\in\Z$ then $\ind_+(P)=\ind(P_k)$ holds independently of $k$.
\end{lemma}
\begin{proof}
Let $k\in\N$ first. Because the first $k$ rows and columns of $P_0$ contain only finitely many nonzero entries, we have the following equality modulo finite rank operators:
\begin{equation} \label{eq:identP+}
P_0\ \cong\ \left(\begin{array}{c|c}0_{k\times k}&0\\\hline0&P_k\end{array}\right)
\ \cong\ \left(\begin{array}{c|c}I_{k\times k}&0\\\hline0&P_k\end{array}\right)
\end{equation}
Consequently, $\ind(P_0)=\ind(P_k)$. For $k\in\Z\setminus\N$ the argument is very similar.
\end{proof}

We will write $P_+$ for any of these singly infinite submatrices $P_k$. Notice that $P_{+}$ can be the zero matrix (not Fredholm) when $P$ is not banded. An example is the permutation that exchanges every pair $i$ and $-i$, for $i\in\Z$.

Here is a concrete formula for the plus-index of permutation matrices that ``split'':
\begin{theorem} \label{thm:split}
Let $\pi$ be a permutation of the integers (not necessarily banded) and denote the corresponding matrix by $P$. Suppose there exist $i^{*}$ and $j^{*}$ such that
\begin{equation} \label{eq:split}
\{\pi(i)\ :\ i<i^{*}\}\ =\ \{j\ :\ j<j^{*}\}.
\end{equation}
Then $\ind_{+}(P)=j^{*}-i^{*}$.
\end{theorem}
\begin{proof}
If such integers $i^{*}$ and $j^{*}$ exist then $P$ decouples into two blocks:
\begin{equation} \label{eq:decouple}
P\ =\  \left(\begin{array}{cc|cc}\ddots&\vdots&\\
\cdots&p_{i^{*}-1,j^{*}-1}&\\\hline
&&p_{i^{*},j^{*}}&\cdots\\
&&\vdots&\ddots
\end{array}\right)
\ =:\ \left(\begin{array}{c|c}P^{(1)}&0\\\hline 0&P^{(2)}\end{array}\right)
\end{equation}
The meeting point $(i^{*}, j^{*})$ may not fall on the zeroth diagonal of $P$ (in fact, it falls on the main diagonal). Because $P$ is invertible, $P^{(1)}$ and $P^{(2)}$ are invertible. It is now easy to show that $j^{*}-i^{*}$ is the plus-index of $P$:

{\bf Case 1.} If $i^{*}\le j^{*}$ then, putting $k:=i^{*}$, we have that $P_{k}$ starts with $j^{*}-i^{*}$ zero columns. Those are followed by $P^{(2)}$, so that the index of $P_{k}$ is $(j^{*}-i^{*})-0=j^*-i^*$.

{\bf Case 2.} If $i^{*}> j^{*}$ then, putting $k:=j^{*}$, we have that $P_{k}$ has $i^{*}-j^{*}$ zero rows, followed by $P^{(2)}$. Again the index is $0-(i^{*}-j^{*})=j^*-i^*$.

\noindent It remains to apply Lemma \ref{lem:Ak}.
\end{proof}

Adding $\kappa=j^{*}-i^{*}$ to all row numbers moves the meeting point of the blocks in \eqref{eq:decouple} to position $(j^{*},j^{*})$. In fact $\ind_{+}(A)=j^{*}-i^{*}$ holds for all invertible matrices $A$ (not just permutations) that decouple in the sense of \eqref{eq:decouple}.

The first question is whether or not such a splitting will appear in every permutation of $\Z$ -- and the unfortunate answer is no!
\begin{example} {\bf \!(symmetric and intertwined) } \label{ex:ex3}
Let us depict this permutation by a graph with vertex set $\Z$ and with a directed edge (an arrow) from $i$ to $j$ iff $j=\pi(i)$.
%
\begin{center}
\includegraphics[width=0.95\textwidth]{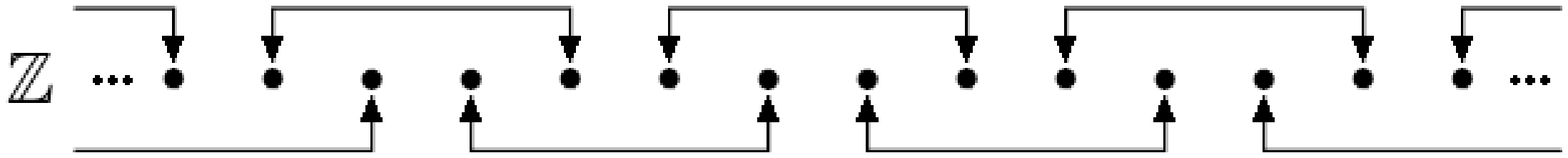}
\end{center}
%
This permutation is symmetric $(\pi(i)=j$ iff $\pi(j)=i)$. In such a case, splitting happens iff the graph falls into separate components $\{i<i^{*}\}$ and $\{i\ge i^{*}\}$. Then the graph can be split in two without cutting any edge, which is obviously impossible in our example.

Of course, this symmetric sitation means that $P=P^{\top}$ and hence $P_{+}=(P_{+})^{\top}$ so that $\ind(P_{+})$ must be zero -- and we don't need help from Theorem \ref{thm:split}. But here is a non-symmetric version of essentially the same example:
\end{example}
\newpage 

\begin{example} {\bf \!(previous example, shifted)} \label{ex:ex4}  $\pi'(i):=\pi(i)+1$.
\begin{center}
\includegraphics[width=0.95\textwidth]{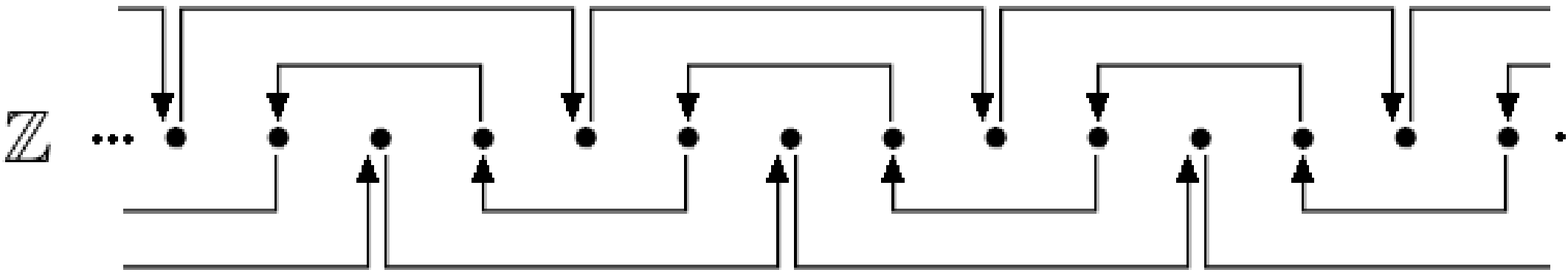}
\end{center}
Without symmetry it is not so easy to see that splitting is impossible. In fact, this graph is disconnected (it has three components) but it cannot be split.
\end{example}

There is a simple trick to make splitting immediately visible in every permutation's graph. We will demonstrate this trick for a simpler example \ref{ex:shift}.

\begin{example} \label{ex:shift}
Let $\pi(i):=i+1$. The corresponding matrix is our shift, $S$.
\begin{center}
\includegraphics[width=0.95\textwidth]{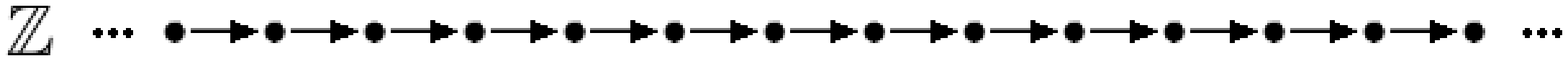}
\end{center}
The graph is clearly connected but we have a split at every single position $i^{*}$, with $j^{*}=i^{*}+1$ in \eqref{eq:split}. This is more easily seen by depicting the $i$'s and the $j$'s separately, where we draw an arrow from $i$ to $j$ if $j=\pi(i)$:
%
%
%
\begin{center}
\includegraphics[width=0.95\textwidth]{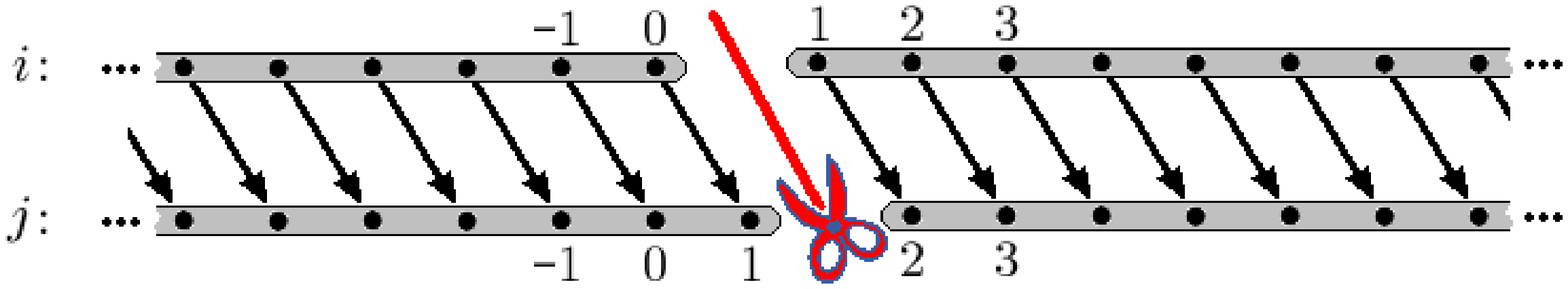}\end{center}
\end{example}

In general, splitting \eqref{eq:split} is equivalent to a separation like this:
\begin{center}
\includegraphics[width=0.95\textwidth]{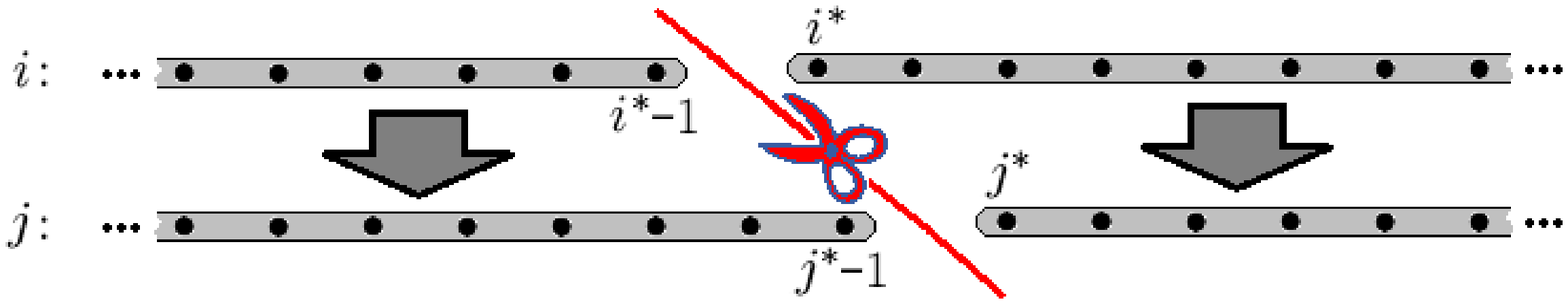}
\end{center}
The nodes $\{i<i^{*}\}$ are only connected to $\{j<j^{*}\}$ and the nodes $\{i\ge i^{*}\}$ only to $\{j\ge j^{*}\}$.
By shifting the $i$-axis accordingly, one can straighten the cut to perfectly vertical. This is exactly what index cancellation does. Shifting the $i$-axis corresponds to renumbering the rows of our matrix. The plus index $j^{*}-i^{*} $ is obvious if one thinks in terms of these $i$-$j$-graphs rather than the matrix.
\noindent The conclusion that there is no split in Example \ref{ex:ex4} is now a simple exercise.
\newpage 

%
%

How do we compute the plus-index when the graph does not split? Here is our agenda:

\begin{itemize}
\item By changing finitely many arrows, we construct a new permutation $\pi'$ that does split.
\item This is possible with the split at any given position $j^{*}$.
\item Changing finitely many arrows only changes finitely many matrix entries. Therefore $\pi$ and $\pi'$ have the same plus-index, and we know that the index for $\pi'$ is $j^{*}-i^{*}$.
\end{itemize}

{\bf Step 1: Delete arrows. } Choose an arbitrary position $j^{*}\in\Z$ where the split should cross the $j$-axis. First delete the $2w$ arrows $i\mapsto\pi(i)$ with $i=j^{*}-w,...,j^{*}+w-1$ from the graph of $\pi$. Then the remaining diagram will have a big gap:
\begin{center}
\includegraphics[width=0.95\textwidth]{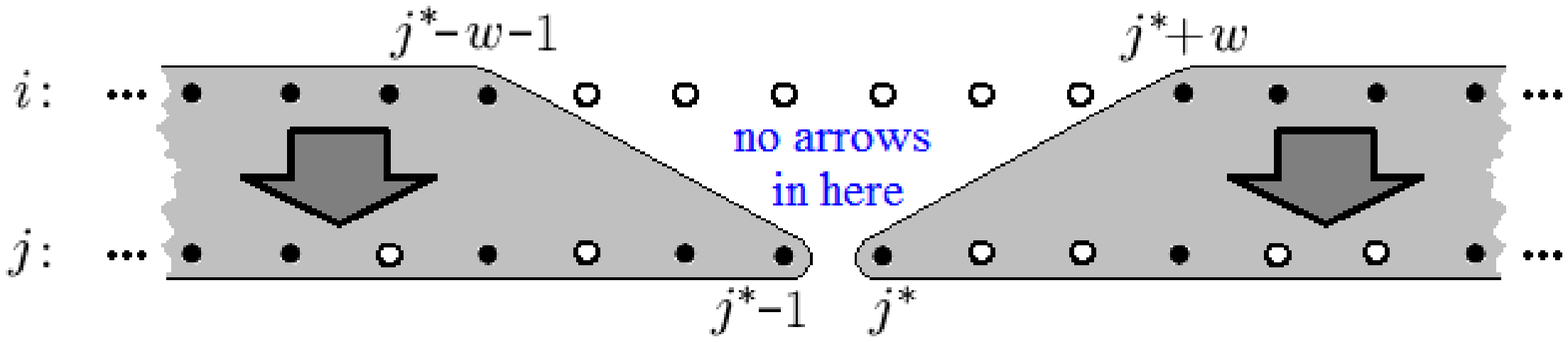}
\end{center}
The endpoints of the $2w$ deleted arrows are now empty. Put those endpoints (the column numbers) into {\sl ascending order}. Suppose $n$ of those column numbers are greater than or equal to $j^*$ -- their endpoints are to the right of the split. The other $2w-n$ endpoints are to the left.

{\bf Step 2: Rewire. } Insert $2w$ new arrows connecting the empty starting points $i=j^{*}-w,...,j^{*}+w-1$ to the ({\sl ordered!}) empty endpoints.
\begin{center}
\includegraphics[width=0.95\textwidth]{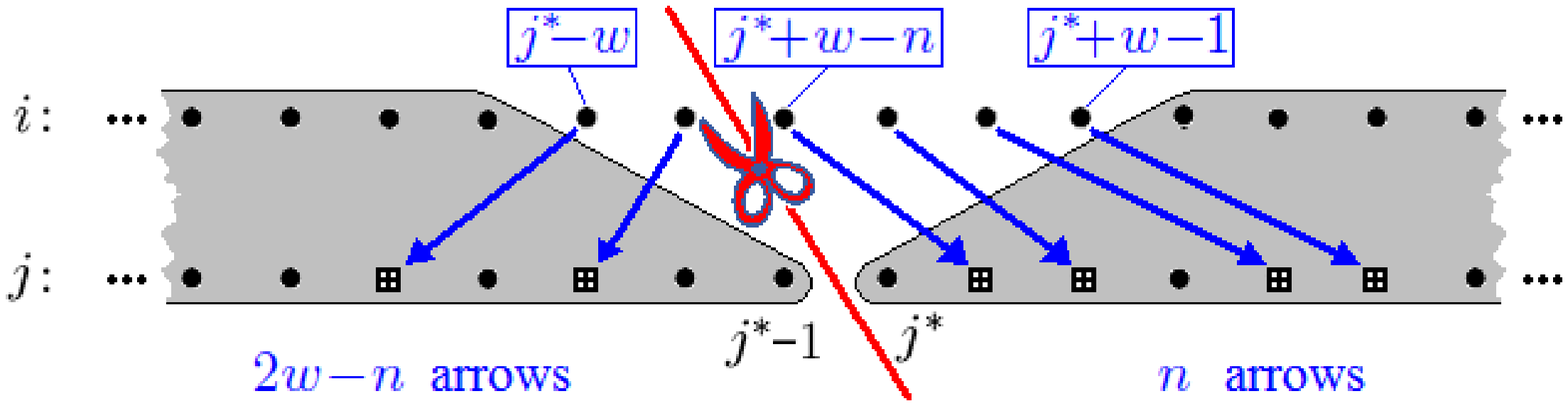}
\end{center}
The new permutation $\pi^{\prime}$ splits at $i^*=j^*+w-n$. If $P^{\prime}$ is the matrix for $\pi^{\prime}$, then $P^{\prime}_+-P_+$ is of finite rank. Therefore $P_+$ and $P^{\prime}_+$ have the same index. This index comes directly from the splitting point of $\pi^{\prime}$:
\begin{equation} \label{eq:indfinal}
\kappa\ =\ \ind_+(P)\ =\ \ind(P_+)\ =\ \ind(P^{\prime}_+)\ =\ j^{*}-i^{*}\ =\ n-w.
\end{equation}
\begin{remark}
We are just reordering the $1$'s in $2w$ consecutive rows of $P$. In the new order, those $1$'s go from left to right. The figure has $w=3$ and there are $n=4$ ones on the right. Then in this example the index is $\kappa=4-3=1$.
\end{remark}

{\sl The number of $1$'s in columns} $j\geq j^{*}$ {\sl is} $n$. Thus we are counting the $1$'s in the following submatrix $R$, and subtracting the bandwidth $w$ to obtain the index $\kappa=n-w$:
\begin{equation} \label{eq:submatrix}
R=\left(\begin{array}{ccc}p_{j^{*}-w,j^{*}}&\cdots&0\\
\vdots& \ddots &\vdots\\
p_{j^{*}+w-1,j^{*}}&\cdots&p_{j^{*}+w-1,j^{*}+2w-1}\end{array}\right)\end{equation}
This submatrix $R$ is lower triangular because $P$ has bandwidth $w$.
This completes our first proof of Theorem $1$.

%

To some extent this result is surprising. The Fredholm index is robust under perturbations in finitely many matrix entries -- so it must be encoded somewhere deep down at infinity. (This is exactly where the limit operator approach \cite{RaRoSiBook,Li:Book,CWLi2008:Memoir}, that we mentioned before, enters the stage.) But for banded permutation matrices, the plus-index can be computed from a finite submatrix -- \textsl{independent of the position}. This is clearly a consequence of the `stiff' rules (exactly one $1$ in each row and column -- and never leave the band!) that these matrices have to follow.


{\bf Second proof of Theorem 1.} The first proof identified the ``right half'' of any $2w$ consecutive rows of $P$, as sufficient to determine the index of $P_+$. Now we can give a direct proof by counting the $1$'s in that right half $R$. In the figure below, the infinite submatrix $P_+$ is marked out by double lines, starting on the zeroth diagonal (marked with $\times$'s). We divide $R$ into a part $R_1$ outside $P_+$ and a part $R_2$ inside $P_+$. Now count the $1$'s in each part\,:
\[
\begin{array}{ccccl@{}l@{\hspace{.2ex}}|ccccccc}
& & & & & & & & & & & & \\[-2.75ex] \cline{2-12}
& & & & & & \bullet & & & & & & \\
{w \textrm{ rows}} & & & & & & \bullet & \bullet & & & & & {R_1}\\
& & & & & & \bullet & \bullet & \bullet & & & & \\[-.2ex] \cline{6-12}
& & & & \vline & & & & & & & & \\[-2.5ex] \cline{2-12}
& & & & \vline & & \times  & \bullet & \bullet & \bullet & & & \\
{w \textrm{ rows}} & & & & \vline & &  \bullet & \times       & \bullet & \bullet & \bullet & & {R_2}\\
& & & & \vline & & \bullet & \bullet & \times       & \bullet & \bullet & \bullet & \\ \cline{2-12}
& & & & \vline & & & & & & & & \\
\end{array}
\]
Each $1$ in $R_1$ will mean a zero column below it in $P_+$. The other columns of $P_+$ are complete:
\[
\textrm{dimension of nullspace of $P_+ =$ number of $1$'s in $R_1$}.
\]
Any zero rows in $P_+$ will be zero rows in $R_2$. Those are counted by the $1$'s that are missing from the $w$ rows of $R_2$:
\[
\textrm{codimension of range of $P_+=w\,-$ (number of $1$'s in $R_2$)}.
\]
Subtraction gives the index formula in Theorem $1$:
\begin{equation}\label{eq:subtraction}
\textrm{index of $P_+ =$ (number of $1$'s in $R$) $-\,w=n-w.$}\quad \proofend
\end{equation}
%

\section{Factorizations of banded permutations}
This short section establishes the factorizations $P_{c}=BC$ and $P_c=F_{1}F_{2}\cdots F_{N}$ of centered permutation matrices $P_c$ (plus-index equal to $0$). The block-diagonal permutations $B$, $C$ have block size $2w$ and $F_{1},\cdots,F_{N}$ have block size $1$ or $2$.

Both factorizations are particularly simple cases of a more general theorem \cite{Strang:GroupsBandInv,Strang:LPU} for banded matrices with banded inverses. For permutations, $P$ and $P^{-1}=P^\top$ have the same bandwidth $w$.
Figure \ref{fig:P0=BC} below shows a typical matrix $P_c$ with $w=2$. Each non-empty square of size $2w=4$ contains two $1$'s (by Theorem 1), because $P_c$ is centered ($\kappa=0$).

{\bf Step 1.} There is a $1$ in each row and column. Row exchanges from the block $B_1^{-1}$ will move the $1$'s in the first square into the top two rows (and thus $1$'s go into the last two rows of the second square). The block $B_2^{-1}$ acts in the same way on the next four rows of $P_c$. All the circled entries $\obullet\,$ have become zero.

{\bf Step 2.} $C_1^{-1}$ executes column exchanges to move the $2+2=4$ ones in its columns to the main diagonal (marked by $\times$). When all blocks $B_i^{-1}$ and $C_i^{-1}$ act in this way, we reach the identity $B^{-1}P_cC^{-1}=I$ and hence $P_c=BC$. This is Theorem $2$.
\begin{figure}[h]
\[
\begin{array}{r|cccc|cccc|cccc|c}
&\multicolumn{4}{|c|}{C_0^{-1}}&\multicolumn{4}{|c|}{C_1^{-1}}&\multicolumn{4}{|c|}{C_2^{-1}}&\\ \hline
& \bullet&\bullet&\times&\bullet & \obullet&&& &&& &\\
& &\bullet&\bullet&\times &\obullet&\obullet&& &&& &\\
\textrm{multiply by }B_1^{-1}& &&\obullet&\obullet & \times&\bullet&\bullet& &&& &\\
& &&&\obullet & \bullet&\times&\bullet&\bullet &&& &\\ \hline
& &&& & \bullet&\bullet&\times&\bullet & \obullet&&& &\\
& &&& & &\bullet&\bullet&\times & \obullet&\obullet&& &\\
\textrm{multiply by }B_2^{-1}& &&& & &&\obullet&\obullet & \times&\bullet&\bullet& &\\
& &&& & &&&\obullet & \bullet&\times&\bullet&\bullet & \\ \hline
& &&& & &&& & &&& &
\end{array}
\]
\caption{{\footnotesize Each square has rank 2. $B_1^{-1}$ and $B_2^{-1}$ produce zeros in the circled positions by row exchanges. Then $C_1^{-1}$ produces (by column exchanges) the 4 by 4 identity matrix with ones in the diagonal positions $\times$.}}
\label{fig:P0=BC}
\end{figure}

A similar factorization was established in \cite{Strang:Fast} for any banded matrix with banded inverse. The same result was found for unitary banded matrices in \cite{RaRoRoe}, with a different proof. It appears that $2w$ is the right block size for these block-diagonal factorizations $A=BC$.

Theorem 3 says that it is possible to factor $P_c$ into even simpler block-diagonal permutations $F_1F_2\cdots F_N$. Now the blocks of each factor are $2\times 2$ or $1\times 1$. Thus each $F_i$ exchanges a set of disjoint pairs of neighbors (as in the bubblesort algorithm, but with disjoint exchanges in parallel for much greater efficiency).

A beautiful proof that $N<2w$ was given by Greta Panova \cite{Panova,Alb,Strang:LPU}. We won't repeat the details. Her key idea is a variation on the ``wiring diagram'' of a permutation.


Briefly, that diagram connects each point $(0,i)$ for $i\in\Z$ by a straight line to the point $(1, \pi_c(i))$. The intersections of those lines tell us the order in which to exchange neighbors. {\sl Key point}: All intersections lie on $N<2w$ vertical lines in the  \cite{Panova,Strang:LPU}. Exchanges on each vertical line can be executed in parallel by a matrix $F_i$ (its bandwidth is $1$). Then $P_c$ is the product $F_1F_2\cdots F_N$ with $N<2w$. 

\part{Band-dominated matrices} \label{part:2}
\section{The plus-index of a band-dominated matrix}
$A$ is a {\sl band matrix} if it is supported on finitely many diagonals only, and it is a {\sl band-dominated matrix} if it is the limit (in the $\ell^{2}\to\ell^{2}$ operator norm) of a sequence of band matrices.

$A_k$ is the semiinfinite submatrix \eqref{eq:defAk} consisting of rows $i\ge k$ and columns $j\ge k$. Then Lemma \ref{lem:Ak} generalizes from permutation matrices to band-dominated matrices that are Fredholm:

\begin{lemma} \label{lem:Ak2}
If $A$ is band-dominated and Fredholm then the index of $A_k$ is independent of $k\in\Z$.
\end{lemma}
\begin{proof}
If $A$ is a band matrix then the proof is literally that of Lemma \ref{lem:Ak} above. If $A$ is band-dominated (i.e.~in the closure of the set of band matrices) then the identification \eqref{eq:identP+} still holds modulo compact operators, which is enough to prove the claim. \end{proof}

In analogy to $A_+=A_1$ we introduce the semiinfinite matrix $A_{-}$ that ends at $a_{0,0}$:
\[
A_-\ :=\ (a_{ij})_{i,j=-\infty}^0\ =\
\left(\begin{array}{ccc}
\ddots&\vdots&\vdots\\
\cdots&a_{-1,-1}&a_{-1,0}\\
\cdots&a_{0,-1}&a_{0,0}
\end{array}\right)
\]
Then $A_{-}$ and $A_{+}$ will together determine the index of $A$:
\begin{equation} \label{eq:A+-}
A\ =\ \left(\begin{array}{c|c}A_-&A_{-+}\\\hline A_{+-}&A_+\end{array}\right)
\ \cong\ \left(\begin{array}{c|c}A_-&0\\\hline 0&A_+\end{array}\right)
\end{equation}
modulo compact operators since $A_{-+}$ and $A_{+-}$ are compact. (Those off-diagonal blocks have finite rank if $A$ is a band matrix and they are norm limits of finite rank matrices, hence compact, if $A$ is band-dominated.)
%
%
From this identification \eqref{eq:A+-} we immediately have Lemma \ref{lem:A+-}:

\begin{lemma} \cite{RaRoRoe} \label{lem:A+-}
A band-dominated biinfinite matrix $A$ is Fredholm iff both $A_+$ and $A_-$ are Fredholm. In that case
\[
\ind(A)\ =\ \ind(A_+)\ +\ \ind(A_-).
\]
\end{lemma}
So if $A$ is invertible then $\ind(A_+)=-\ind(A_-)$. Clearly, all results that we prove for the plus-index have their counterpart for the minus-index, $\ind_-(A)=\ind(A_-)$. Here is another important result on the plus-index:
\begin{lemma} \label{lem:AB+}
If $A$ and $B$ are band-dominated and Fredholm biinfinite matrices then
\[
\ind_+(AB)\ =\ \ind_+(A)\ +\ \ind_+(B).
\]
\end{lemma}
\begin{proof}
By \eqref{eq:A+-} we have, modulo compact operators,
\[
AB\ \cong\ \left(\begin{array}{c|c}A_-&0\\\hline 0&A_+\end{array}\right)
\left(\begin{array}{c|c}B_-&0\\\hline 0&B_+\end{array}\right)\ =\
\left(\begin{array}{c|c}A_-B_-&0\\\hline 0&A_+B_+\end{array}\right),
\]
so that $(AB)_+\cong A_+B_+$, whence $\ind((AB)_+)=\ind(A_+B_+)=\ind(A_+)+\ind(B_+)$.
\end{proof}

Recalling the forward shift $S$ with $\ind_+(S)=0-1=-1$, Lemma \ref{lem:AB+} yields
\begin{equation} \label{eq:IC}
\ind_+(S^\kappa A)\ =\ \kappa\cdot \ind_+(S)\ +\ \ind_+(A)\ =\ -\kappa\ +\ \ind_+(A).
\end{equation}
Then if $\kappa=\ind_+(A)$, index cancellation $A_c:=S^\kappa A$ indeed leads to
plus-index zero.
\section{Triangular matrices}
The plus-index of a general band-dominated matrix $A$ is not easy to compute. We will assume that $A$ is banded and invertible (as a bounded operator on $\ell^{2}(\Z)$). One possible approach (far from complete in this paper) is to factor $A$ into a lower triangular $L$ and an upper triangular $U$, with a banded permutation $P$ in between:
\[
A=LPU=(\hbox{lower triangular}) (\hbox{permutation}) (\hbox{upper triangular}).
\]
\noindent This extends the factorization that comes from Gaussian elimination on a finite invertible matrix. Notice the ``Bruhat convention'' that places $P$ between $L$ and $U$. In that position $P$ is unique. When $P$ is a finite matrix, the $1$'s are all determined by the ranks of the upper left submatrices of $A$. Elimination proceeds as normal (subtracting multiples of the pivot rows from lower rows) except that we wait to the end to reorder the pivot rows into $U$ by using $P$.

The steps are described in \cite{Strang:TriangFact}, where the main purpose is to extend $A=LPU$ to biinfinite matrices $A$ (banded and invertible). The factors $L,P,U$ have bandwidth $\leq 2w$. But the elimination process has to be reconsidered,  since it can no longer start at $a_{11}$ (which is not the first entry, it is in the center of the matrix). Briefly, we look at the upper left submatrices $A_{k-}$ (singly infinite) containing columns $\leq k$ and rows $\leq k+w$. By Fredholm theory, the index of $A_{k-}$ is independent of $k$. Its columns are independent because $A$ is invertible (and all nonzeros survive into $A_{k-}$). Any dependent rows of $A_{k-}$ are among the last $2w$ rows (for the same reason). So the minus-index of $A$ and the location of its main diagonal are determined by the number $d$ of dependent rows.

Since $d$ does not depend on $k$, one row changes to independent and a new dependent row appears (if $d>0$) when $k$ increases to $k+1$. That newly independent row is the pivot row in the following elimination step. Dependence involves {\sl all} the earlier rows of $A_{k-}$, so this elimination process is not constructive---at least not in the usual sense.

{\sl Note:} We say that an infinite set of vectors $\{v_i\}_{i\in\I}$ with $\I\subset\Z$ is {\sl linearly dependent} if there is a non-zero sequence $(c_i)_{i\in\I}\in\ell^2(\I)$ such that $\sum_{i\in\I} c_iv_i=0$. In that sense, the columns of an infinite matrix form a linearly independent set iff that matrix is injective as an operator on $\ell^2$.

Returning for a moment to permutation matrices, rows of $P_{k-}$ are dependent when they are zero. They are independent when they contain a $1$. Therefore the number $d$ of zero rows (in the left half of any $2w$ rows of $P$) involves the same count of $1$'s as in Theorem $1$ (in the right half of those rows).
\newpage 

Now suppose that $A$ is invertible and banded, with $A=LPU$, where all three factors have bounded inverses. (This can fail even for block-diagonal matrices $A$ with $2\times 2$ orthogonal blocks. The upper left entries in those blocks can approach zero.) As usual, $L$ and $L^{-1}$ are lower triangular, $P$ is a permutation, and $U$ and $U^{-1}$ are upper triangular.  Then $P$ contains all information about indices and the correct position of the main diagonal.

\begin{lemma} \label{lem:LPU}
The plus-index of $A$ equals the plus-index of $P$ (and that is easily computable from $P$). Also the minus-indices are equal.
\end{lemma}
\begin{proof}
As shown in \cite{Strang:TriangFact}, $L$ and $U$ are both banded, with bandwidth $\le 2w$.
By Corollary \ref{cor:triangular} below, $\ind_+(L)$ and $\ind_+(U)$ are zero. It follows from Lemma \ref{lem:A+-} that also $\ind_-(L)$ and $\ind_-(U)$ are zero.
Now use the key property that $L$ and $U$ are triangular:
\[
\left(
\begin{array}{@{\;}ll@{\;}}
	A_{-} & A_{-+}\\
	A_{+-} & A_{+}\end{array}
\right)=\left(\begin{array}{@{\;}ll@{\;}}
	L_{-} & 0\\
	L_{+-} & L_{+}\end{array}\right)
\left(
\begin{array}{@{\;}ll@{\;}}
P_{-} & P_{-+}\\
P_{+-} & P_{+}
\end{array}\right)\left(
\begin{array}{@{\;}ll@{\;}}
	U_{-} & U_{-+}\\
	0 & U_{+}
\end{array}\right).
\]
The upper left block $A_{-}$ immediately factors into
\begin{equation}
	A_{-}=L_{-}P_{-}U_{-}.
	\label{eq:A+}
\end{equation}
%
By the formula for the index of a product we get
\begin{equation} \label{eq:A+_}
\ind(A_{-})\ =\ \ind(L_-)+\ind(P_{-})+\ind(U_-)\ =\ \ind(P_-)
\end{equation}
and then $\ind_+(A)=-\ind_-(A)=-\ind_-(P)=\ind_+(P)$, again by Lemma \ref{lem:A+-} since $A$ and $P$ have index zero (they are invertible).
\end{proof}

Alternatively, one can apply Lemma \ref{lem:AB+} directly to $A=LPU$ to get
\[
\ind_+(A)\ =\ \ind_+(LPU)\ =\ \ind_+(L)+\ind_+(P)+\ind_+(U)\ =\ \ind_+(P)
\]
and then conclude $\ind_-(A)=\ind_-(P)$ from Lemma \ref{lem:A+-}.
\medskip

We continue with some results on lower triangular biinfinite matrices and their plus-index (which locates their main diagonal). It is not surprising, but however in need of a proof, that the main diagonal of a lower triangular matrix has to be in the lower triangle. This is what we show now. So let $A$ be a biinfinite matrix.

\begin{lemma} \label{lem:B+}
If $A$ is lower triangular and Fredholm (but not necessarily band-dominated) then, for sufficiently large $k$, the submatrix $A_k$ from \eqref{eq:defAk} is injective (as an operator on $\ell^2(\N)$).
\end{lemma}
\noindent{\bf First proof.}
If $A$ is Fredholm and $\alpha=\dim(N(A))$ then there is a set $J=\{j_1,...,j_\alpha\}$ of integers such that the columns of $A$ whose number is not in $J$ form a linearly independent set. Let $k$ be bigger than $max(J)$. Then all columns $j\ge k$ of $A$ and hence (because $A$ is lower triangular) all columns of the semiinfinite submatrix $A_k$ are linearly independent. So $A_k$ is injective.
\proofend

\noindent{\bf Second proof. }
Since $A$ is Fredholm on $\ell^2$, there are operators $B$ and $K$ such that $BA=I+K$ with $B$ bounded and $K$ compact on $\ell^2$ \cite{Davies2007:Book,GohGoldKash}. For $k\in\Z$, let $I_k$ denote the operator $\ell^2\to\ell^2$ that puts all entries $x_i$ of $x\in\ell^2$ with $i<k$ to zero and leaves all $x_i$ with $i\ge k$ unchanged. ($I_k$ is \eqref{eq:defAk} for $A=I$.) Since $A$ is lower triangular, we have $I_kAI_k=AI_k$ for all $k\in\Z$. Moreover, $\|I_k x\|_{\ell^2}\to 0$ as $k\to +\infty$ for all $x\in\ell^2$. By compactness of $K$, it follows \cite{GohGoldKash} that the operator norm $\|I_kK\|$ goes to zero as $k\to +\infty$. So fix $k\in\Z$ large enough that $\|I_kK\|<1$. Then $C=I+I_{k}K$ is invertible (by Neumann series). Now
\[
I_kBI_kAI_k\ =\ I_kBAI_k\ =\ I_k(I+K)I_k\ =\ I_k+I_kKI_k\ =\ (I+I_{k}K)I_{k}\ =\ CI_{k},
\]
so that $C^{-1}I_kBI_kAI_k=I_{k}$. Now take $x$ from the range of $I_k$ such that $I_kAI_kx=0$.
Then $0=C^{-1}I_kBI_kAI_kx=I_{k}x=x$, so that $x=0$ is the only $x$ in the range of $I_k$ with $I_kAI_kx=0$. By
$I_{k}AI_{k}|_{R(I_{k})}=A_{k}$ this means that the only solution $y=(y_i)_{i=k}^\infty$ of $A_ky=0$ is the trivial solution $y=0$.
\proofend

\begin{corollary} \label{cor:triangular1}
If $A$ is band-dominated, lower triangular and Fredholm then $\ind_+(A) \le 0$.
\end{corollary}
\begin{proof}
By Lemma \ref{lem:B+}, the operator behind $A_k$ has nullspace $\{0\}$ for all sufficiently large $k\in\Z$, so that $\ind(A_k)=0-\beta\le 0$. The claim now follows from Lemma \ref{lem:Ak2}.
\end{proof}

So indeed, the main diagonal of a lower triangular matrix is in the lower triangle.
By simple translation via $S^{d}$, the main diagonal of a matrix that is zero above the $d$-th diagonal must be on or below that $d$-th diagonal. By passing to the adjoint matrix, one can write down an analogous statement for upper triangular matrices and their translates and then combine the two: The main diagonal of a band matrix must be in that band.

Consequently, tridiagonal matrices have plus-index $\kappa=-1,\,0$ or $1$. Using this fact and our Theorem 1, it is easy to see that the only tridiagonal permutation matrices that are not centered (i.e.~they have a nonzero plus-index) are the shifts $S$ and $S^{-1}$. Indeed, suppose $P$ is a biinfinite permutation matrix with bandwidth $w=1$ and plus-index $\kappa=-1$. Then,
for every $j^*\in\Z$, the number of ones in the $2\times 2$ matrix \eqref{eq:submatrix} is $n=\kappa+w=0$. So the $1$ in row $j^*$ must be in column $j^*-1$, whence $P$ must be $S$.
Similarly, $w=1$ and $\kappa=1$ leads to $P=S^{-1}$.

Here is another consequence of Corollary \ref{cor:triangular1}:
\begin{corollary} \label{cor:triangular}
If $A$ is band-dominated, lower triangular and invertible with a lower triangular inverse then $\ind_+(A)=0$.
\end{corollary}
\begin{proof}
By Corollary \ref{cor:triangular1}, we have $\kappa=\ind_+(A)\le 0$ and $\lambda=\ind_+(A^{-1})\le 0$. But then $\kappa=-\lambda$ by Lemma \ref{lem:AB+}, so that $\kappa$ must be zero.
\end{proof}

The latter shows how Problem $1$ from the beginning of our paper is related to the plus-index.
After discussing Problems $1$ and $2$, we mentioned that they are related: If the finite section method applies to $A$ then a lower triangular $A$ will have a lower triangular inverse.
Now we discuss an amazing extension of this statement, coming from the following theorem (for finite matrices) by Asplund \cite{Asplund}\,:
\newpage 

\begin{theorem} \label{thm:Asplund}
Let $A$ be an invertible matrix and fix two integers $p$ and $k$. Then the following are equivalent:\\[-4ex]
\begin{itemize}
\item[(i)] All submatrices $B$ above the $p$-th superdiagonal of $A$ have $\rank(B)<k$.\\[-4ex]
\item[(ii)] All submatrices $C$ above the $p$-th subdiagonal of $A^{-1}$ have $\rank(C)<p+k$.
\end{itemize}
\end{theorem}
See \cite{Strang:Asplund} for discussion and a new proof. With $p=0$ and $k=1$ we get the familiar statement that lower triangular matrices have lower triangular inverses. With $p=1$ and $k=1$, one can see that all submatrices above the first subdiagonal (and similarly: below the first superdiagonal) of $A^{-1}$ have rank $<2$ if $A$ is tridiagonal.

Let us attempt to transfer Asplund's theorem to singly and doubly infinite matrices. Our assumption will be that the {\sl finite section method} (short: {\sl FSM}) applies to the infinite matrix $A$. Here is again what that means:

The {\sl finite sections} of $A$ are square submatrices $A_{n}$ whose upper left and lower right corners lie at positions $l_n$ and $r_n$ on the zeroth diagonal of $A$. The sequences $l_n$ and $r_n$ go to $-\infty$ and $+\infty$ respectively, except when $A$ is only semi-infinite -- then $l_n$ is fixed at $1$ and $r_n$ goes to $+\infty$. One says that the FSM applies to $A$ if: $A$ is invertible, the matrices $A_n$ are invertible for sufficiently large $n$, and their inverses converge strongly to $A^{-1}$.
This implies that $A$ is centered (see \cite{Li:FSMsubs,SeidelSilbermann3}). Also note that strong convergence ($A_n^{-1}x\to A^{-1}x$ for all $x$) implies entrywise convergence of the matrices $A_n^{-1}$ to $A^{-1}$. See e.g.~\cite{Li:Book,Li:FSMsubs,RaRoSiBook,Roch:FSM,SeidelSilbermann3} and the references therein for more on the FSM.

Under this (reasonable) assumption, we now extend Asplund's Theorem \ref{thm:Asplund} to a singly or doubly infinite matrix $A$ (not necessarily band-dominated):

\begin{proof}
Suppose $(i)$ holds, i.e.~all (finite and infinite) submatrices $B$ above the $p$-th superdiagonal of $A$ have $\rank(B)<k$. To show that $(ii)$ holds, let $C$ be an arbitrary (finite or infinite) submatrix above the $p$-th subdiagonal of $A^{-1}$. We show that $\rank(C)<p+k$.

If $C$ is infinite then $\rank(C)$ is the supremum of all $\rank(C')$ with $C'$ going through all {\sl finite} submatrices of $C$. (In particular, $\rank(C)=\infty$ iff the set of all $\rank(C')$ is unbounded.) Therefore it is enough to show that $\rank(C')<p+k$ for all {\sl finite} submatrices of $A^{-1}$ above the $p$-th subdiagonal. So we can assume that $C$ is a finite matrix.

Let $l_1,l_2,...$ and $r_1,r_2,...$ be the cut-off positions for the (by our assumption applicable) finite sections of $A$. For all sufficiently large $n$ (say $n>N$) the interval $[l_n,r_n]$ contains the row and column numbers in which $C$ is positioned at $A^{-1}$. By applicability of the FSM, the inverses of our finite sections $A_n$ of $A$ converge entrywise to $A^{-1}$. Let $C_n$ denote the submatrix of $A_n^{-1}$ that is at the same position as $C$ is in $A^{-1}$, so that $C_n\to C$ entrywise as $n\to\infty$.

Now we apply Theorem \ref{thm:Asplund} to the {\sl finite} matrix $A_n$: All submatrices above the $p$-th superdiagonal of $A_n$ have rank $<k$, by $(i)$, because they are submatrices of $A$ (above the $p$-th superdiagonal). So $C_n$, being a submatrix above the $p$-th subdiagonal of $A_n^{-1}$, has rank $<p+k$. And this is true for all $n>N$. For finite matrices, entrywise convergence is convergence in all matrix norms, so that $\|C_n-C\|\to 0$. Rank is a lower semi-continuous function with respect to the matrix norm, so
\[
\rank(C)\ \le\ \liminf_{n\to\infty} \rank(C_n)\ <\ p+k
\]
and we are done with $(i)\Rightarrow (ii)$. The other direction is checked analogously.
\end{proof}

\begin{remark}
Our proof shows that applicability of the FSM is sufficient for Asplund's theorem to hold for an infinite matrix. To see that it is not necessary, go back to the permutation matrix $P$ in Example \ref{ex:ex3}. Asplund's theorem can be seen to hold for this matrix (note that $P^{-1}=P^\top=P$) but the FSM does not apply: No one of the finite sections is invertible because there are no places $l_n$ and $r_n$ to cut without breaking one of the links in this permutation graph.
\end{remark}

%
%
\bigskip

{\bf Acknowledgements. } We would like to thank Steffen Roch for helpful discussions. Moreover, the first author acknowledges the financial support by the Marie-Curie Grant PERG02-GA-2007-224761 of the EU.

\pagebreak
\bigskip
\noindent {\bf Authors:}\\[3mm]
%
%
\noindent Marko Lindner\hfill \href{mailto:marko.lindner@mathematik.tu-chemnitz.de}{{\tt marko.lindner@mathematik.tu-chemnitz.de}}\\
TU Chemnitz\\
Fakult\"at Mathematik\\
D-09107 Chemnitz\\
GERMANY\\[5mm]
\noindent Gilbert Strang\hfill \href{mailto:gilstrang@gmail.com}{{\tt gilstrang@gmail.com}}\\
MIT\\
Room 2-240\\
Cambridge MA 02139\\
USA
\end{document}